\documentclass[11pt]{article}
\usepackage[all]{xy}
\usepackage{amsmath}
\usepackage{amsfonts}
\usepackage{amssymb}
\usepackage{amscd}
\usepackage{amsthm}
\usepackage{latexsym}
\usepackage{amsbsy}

\newtheorem{example}{Example}[section]

\newtheorem{remark}{Remark}

\usepackage[all]{xy}
\usepackage{amsmath}
\usepackage{amsfonts}
\usepackage{amssymb}
\usepackage{amscd}
\usepackage{amsthm}
\usepackage{latexsym}
\usepackage{amsbsy}
\begin{document}

\title{The Algebra of  Formal Twisted Pseudodifferential Symbols and a Noncommutative Residue }
\author {Farzad Fathi-Zadeh and Masoud Khalkhali  \\ Department of Mathematics,  University of Western Ontario \\ London, Ontario, Canada}
\date{}
\maketitle
\begin{abstract} We extend the Adler-Manin trace on the algebra of
pseudodifferential symbols  to a twisted setting.

\end{abstract}

\section{Introduction}
The definition of the noncommutative residue of Adler-Manin
\cite{adl, man} for  the algebra of formal pseudodifferential
symbols on the circle and its full analytic extension to higher
dimensions by Wodzicki and Guillemin \cite{wod1, gui},  relies on
the existence of an invariant trace on the ground algebra. A vast
generalization of this noncommutative residue has been obtained by Connes and
Moscovici in the context of   spectral triples and the local index
formula  in noncommutative geometry in \cite{conmos0}. Here again
the conditions needed to define the residue trace on the algebra
of pseudodifferential operators for the spectral triple,  imply
the existence of a trace on the base algebra. To deal with  `type
III situations' where there can be no trace on the  ground
algebra, Connes and Moscovici introduced the notion of a twisted
spectral triple in \cite{conmos}. This of course raises the
questions of extending the results of \cite{conmos0} to this new
twisted setting.

In this paper we take a  modest first step in this
direction and look for analogues of the Adler-Manin trace on an
algebra of formal twisted pseudodifferential symbols. To our
surprise, we found that starting with an {\it invariant  twisted trace} on the
ground algebra, one can still induce a trace on the algebra of
twisted symbols. This of course suggests that in the context of
twisted spectral triples  one might be able to use the twisted
trace  defined in \cite{conmos}  to define  a trace
on a twisted pseudodifferential calculus for the spectral triple.
We shall not discuss this problem in this paper any further.

Our main results are as follows. For a triple $(A, \sigma, \delta)$ consisting of an algebra $A$,
an automorphism $\sigma: A \to A$,  and a twisted derivation
$\delta: A \to A$, we define an algebra of {\it formal twisted
pseudodifferential symbols} $\Psi (A, \sigma, \delta)$ with
coefficients in $A$. Our construction extends those of Adler-Manin
as well as the algebra of $q$-pseudodifferential symbols on the
circle due to Khesin-Lyubashenko-Roger \cite{klr} and Kassel
\cite{kas2}. We show that starting with an invariant  twisted trace on $A$,
the noncommutative residue functional is a trace on $\Psi (A,
\sigma, \delta)$. We also define an $n$-dimensional analogue of $\Psi (A,
\sigma, \delta)$
  and show that starting with an invariant $\sigma^n$-trace
on $A$, one obtains a trace on $\Psi (A,
\sigma, \delta)$.  We also extend  the
 logarithmic cocycle of \cite{klr, krakhe} to our twisted set up. One
 can
get  a different kind of result if one assumes that $\delta$ and
$\sigma$ commute. In this case the automorphism $\sigma$ extends
to an automorphism of $\Psi (A, \sigma, \delta)$ and we show that
starting with an invariant trace on $A$, the noncommutative
residue is in fact a twisted trace on $\Psi (A, \sigma, \delta)$. In the last  section of this paper
 we give a general method to construct twisted
derivations, twisted traces, and twisted connections on a crossed
product algebra. We also indicate a relation with a twisted spectral
triple constructed in \cite{conmos}.

We should mention that when  we discussed the idea of our algebra of
twisted pseudodifferential symbols with Nigel Higson at a recent
conference at the Fields Institute in Toronto,
he brought a paper of Schneider and
Venjakob  \cite{schven} to our attention where a similar algebra of  `skew power series rings'
has been defined   in the context of noncommutative Iwasawa theory. The
two algebras are however mildly different as we use a different
completion from \cite{schven}.
 Apart from this, the contents and focus  of the two
papers are quite disjoint. We would like to thank the Hausdorff Research
Institute for Mathematics for their generous support and a very nice
environment during our stay in the summer of  2008 when   most of this work was done. 

\section{Preliminaries}

To motivate our definitions in the next section, we recall some standard notions and constructions in this section. In this paper by an algebra we mean an associative, unital, not necessarily commutative algebra over $\mathbb{C}$. The choice of $\mathbb{C}$ as the ground field is not essential and in fact all of our definitions and results can be phrased over an arbitrary field of characteristic zero. The characteristic zero assumption however would be essential. Let $A$ be an algebra. Recall that a \emph{derivation} on $A$ is a linear map $\delta : A \to A$ such that $\delta (ab)= a \delta (b) + \delta (a) b$ for all $a,b \in A$. Given a pair $(A, \delta)$ as above, the algebra of \emph{formal differential symbols} $D(A, \delta)$ \cite{kas1,man} is, by definition, the algebra generated by $A$ (as a subalgebra) and a symbol $\xi$ subject to the relations

\begin{equation} \label{eq:relation}
\xi a - a \xi = \delta (a)
\end{equation}
for all $a \in A$. Every element of $D(A, \delta)$ has a unique expression of the form
\[
D= \sum_{i=0}^{N} a_i \xi^i
\]
for some non-negative integer $N$ and $a_i \in A$. We think of $D$ as a differential operator of order at most $N$. Using \eqref{eq:relation}, one can inductively prove that
\begin{equation} \label{eq:relation2}
\xi^n a = \sum_{j=0}^{n} \binom{n}{j} \delta^j(a) \xi^{n-j}
\end{equation}
for all $a \in A$ and $n \geq 0$. Using \eqref{eq:relation2}, we obtain the following multiplication formula in $D(A, \delta)$:
\[
\Big ( \sum_{i=0}^{M}a_i \xi^i \Big ) \Big ( \sum_{j=0}^{N} b_j \xi^j \Big ) = \sum_{n=0}^{M+N} \Big ( \sum_{i,j,k}^{} \binom{i}{k} a_i \delta^{k}(b_j) \Big ) \xi^n
\]
where the internal summation is over all $0 \leq k  \leq i \leq M$, and $0 \leq j \leq N$ such that $i+j-k=n$.\\

The algebra of \emph{formal pseudodifferential symbols} of $(A, \delta)$ \cite{adl,man,kas1}, denoted by $\Psi (A,\delta)$, is obtained from $D(A, \delta)$ by formally inverting $\xi$ and completing the resulting algebra. More precisely, it is defined as follows. Elements of $\Psi (A,\delta)$ consist of formal sums
\[
D=\sum_{i=- \infty}^{N} a_i \xi^i
\]
with $a_i \in A$, and $N \in \mathbb{Z}$. Its multiplication is defined by first extending \eqref{eq:relation2} to all $n \in \mathbb{Z}$:

\begin{equation}
\xi^n a = \sum_{j=0}^{\infty} \binom{n}{j} \delta^j(a) \xi^{n-j}. \nonumber
\end{equation}
Here the binomial coefficient $\binom{n}{j}$ for $n \in \mathbb{Z}$, and non-negative $j \in \mathbb{Z}$, is defined by $\binom{n}{j}:= \frac{n(n-1) \cdots (n-j+1)}{j!}$. Notice that for $n < 0$, we have an infinite formal sum. In general, for $D_1 = \sum_{i=- \infty}^{M} a_i \xi^i$, and $D_2 = \sum_{j=- \infty}^{N} b_j \xi^j$, the multiplication is defined by
\[
D_1 D_2 = \sum_{n=-\infty}^{M+N} \Big ( \sum_{i,j,k} \binom{i}{k} a_i \delta^{k}(b_j) \Big ) \xi^n
\]
where the internal summation is over all integers $i \leq M$, $ j \leq N$, and $ k \geq 0 $ such that $i+j-k =n$.

Let now $\tau : A \to \mathbb{C}$ be a $\delta$-\emph{invariant trace} on $A$. Thus, by definition, $\tau$ is a linear map and we have $\tau (ab) = \tau (ba)$ and $\tau (\delta (a))=0$ for all $a,b \in A$. The \emph{Adler-Manin noncommutative residue} \cite{adl,man,kas1},  $\text{res} : \Psi (A, \delta) \to \mathbb{C}$ is defined by
\[
\text{res} \, \big ( \sum_{i=-\infty}^{N}a_i \xi^i \big ) = \tau (a_{-1}).
\]
One checks that $\text{res}$ is indeed a trace, i.e.
\[
\text{res} \, [D_1,D_2] = 0
\]
for all $D_1, D_2 \in \Psi (A, \delta)$. Equivalently, one shows that the map
$\text{res} : \Psi (A, \delta) \to A/([A,A]+ \text{im} \, \delta)$
\[D \mapsto a_{-1} \,\,\,\, \text{mod} \,\,\, [A,A] + \text{im} \, \delta
\]
is a trace on $\Psi(A,\delta)$ with values in $A/([A,A] + \text{im} \, \delta)$. 

A relevant example is when $A=C^{\infty}(S^1)$, the algebra of smooth functions on the circle  $S^1= \mathbb{R} / \mathbb{Z}$, with $\tau (f) = \int_{0}^{1} f(x)dx$ and $\delta(f)=f'$. Then the noncommutative residue coincides with Wodzicki residue on the algebra of classical pseudodifferential symbols. \\

\section{Algebra of formal twisted pseudodifferential symbols}

In this section we extend the algebra of formal pseudodifferential symbols to a twisted setup, and look for traces and twisted traces on these algebras.

\newtheorem{twisteddef}{Definition}[section]
\begin{twisteddef}\label{twisteddef}
Let $\sigma : A \to A$  be an automorphism of an algebra $A$.
\begin{enumerate}
\item A $\sigma$-derivation on $A$ is a linear map $\delta : A \to A$ such that
\begin{equation}
\delta (ab) = \delta (a) b + \sigma (a) \delta(b), \qquad \forall a,b \in A. \nonumber
\end{equation}
\item A $\sigma$-trace on $A$ is a linear map $\tau :A \to \mathbb{C}$ such that
\begin{equation}
\tau (ab) = \tau (\sigma (b)a), \qquad \forall a,b \in A. \nonumber
\end{equation}
\end{enumerate}
\end{twisteddef}

Given a triple $(A, \sigma, \delta)$ as in the above definition, we define \emph{the algebra of formal twisted pseudodifferential symbols}, denoted by $\Psi (A, \sigma, \delta)$, as the set of all formal series in $\xi$ with coefficients in $A$:
\[
\Psi (A, \sigma, \delta) := \Big\{ \sum_{n=- \infty}^N a_n  \xi^{n} \mid N \in \mathbb{Z}, \quad a_n \in A \quad \forall n \leq N\Big\}.
\]
To define the multiplication, we impose the relations

\begin{equation} \label{eq:twistedrelation}
\xi a - \sigma (a) \xi = \delta (a)
\end{equation}
for all $a \in A$, and

\begin{equation}
\xi^{-1} \xi = \xi \xi^{-1}=1. \nonumber
\end{equation}
Using \eqref{eq:twistedrelation}, one can inductively show that (\emph{cf.} also \cite{schven}) for all $n \geq 0$,

\begin{equation} \label{eq:twistedrelation2}
\xi^n a = \sum_{i=0}^{n} P_{i, \, n}(\sigma , \delta)(a) \xi^{i},
\end{equation}
where $P_{i,\, n}(\sigma , \delta) : A \to A$ is the noncommutative polynomial in $\sigma$ and $\delta$ with $\binom{n}{i}$ terms of total degree $n$ such that the degree of $\sigma$ is $i$, e.g.

\[
P_{3,4}(\sigma, \delta) = \delta \sigma^3 + \sigma \delta \sigma^2 + \sigma^2 \delta \sigma + \sigma^3 \delta.
\]
We need an extension of \eqref{eq:twistedrelation2} for $n<0$. For $n=-1$, we have by induction:
\[
 \xi^{-1} a= \sum_{i=0}^{N} (-1)^i \sigma^{-1} (\delta \sigma^{-1})^{i}(a) \xi^{-1-i} + (-1)^{N+1} \xi^{-1} (\delta \sigma^{-1})^{N+1}(a)\xi^{-1-N}.
\]
This suggests putting
\begin{equation} \label{eq:twistedrelation-1}
\xi^{-1} a= \sum_{i=0}^{\infty} (-1)^i \sigma^{-1} (\delta \sigma^{-1})^{i}(a) \xi^{-1-i}.
\end{equation}
Multiplying \eqref{eq:twistedrelation-1} by $\xi^{-1}$, we obtain
\[
\xi^{-2} a = \sum_{i=0}^{\infty} \sum_{j=0}^{\infty} (-1)^{i+j} \sigma^{-1} (\delta \sigma^{-1})^j \sigma^{-1} (\delta \sigma^{-1})^i (a) \xi^{-2-i-j}.
\]
Continuing we have
\begin{equation} \label{eq:twistedrelation-n}
\xi^{-n}a = \sum_{i_1=0}^{\infty} \cdots \sum_{i_{n}=0}^{\infty} (-1)^{i_1  + \cdots + i_{n}} \sigma^{-1} (\delta \sigma^{-1})^{i_{n}} \cdots \sigma^{-1} (\delta \sigma^{-1})^{i_1} (a) \xi^{-n - i_1 - \cdots - i_{n}}.
\end{equation}
for any integer $n \geq 1$.
This suggests the following formula for the multiplication of twisted pseudodifferential symbols $D_1= \sum_{n=- \infty}^{N} a_n \xi^n$, and $D_2 = \sum_{m=- \infty}^{M} b_m \xi^m$;
\begin{eqnarray} \label{twistedmultiplication}
D_1 D_2=
&& \sum_{m= - \infty}^{M} \sum_{n<0} \sum_{i \geq 0} (-1)^{| i |} a_n \sigma^{-1} (\delta \sigma^{-1})^{i_{-n}} \cdots \sigma^{-1} (\delta \sigma^{-1})^{i_1} (b_m) \xi^{m+n- |i|} {}
                         \nonumber \\
&& + \sum_{m= - \infty}^{M} \sum_{n=0}^{N} \sum_{j=0}^{n} a_n P_{j,\, n} (\sigma, \delta) (b_m) \xi^{m+j}.
\end{eqnarray}
where $i=(i_1, \dots , i_{-n})$ is an $n$-tuple of integers and $|i|= i_1 + \cdots + i_{-n}$. One can show that, endowed with the multiplication defined in \eqref{twistedmultiplication}, $\Psi(A, \sigma, \delta)$ is an associative unital algebra.

\newtheorem{nctracelemma}[twisteddef]{Lemma}
\begin{nctracelemma} \label{nctracelemma}

Let $A, \sigma, \delta, \tau$ be as in Definition \ref{twisteddef}. If $\tau \circ \delta = 0$, then for any $a,b \in A$, and any m-tuple $i=(i_1, \dots , i_m)$ of non-negative integers, we have:
\begin{equation} \label{eq:nclemmaeq}
\tau \big ( b \sigma^{-1} (\delta \sigma^{-1})^{i_1} \cdots \sigma^{-1} (\delta \sigma^{-1})^{i_m} (a) \big )
 = (-1)^{i_1 + \cdots + i_m} \tau \big( a \delta^{i_m} \sigma \delta^{i_{m-1}} \cdots \sigma \delta^{i_1} (b) \big)
\end{equation}

\begin{proof}
Since $\tau$ is a $\sigma$-trace, we have
\begin{equation} \label{nclemmaeq0}
\tau(\sigma (a)) = \tau(a), \,\,\, \forall a \in A.
\end{equation}
From $\sigma$-derivation property of $\delta$ and $\tau \circ \delta=0$ , it follows that:
\begin{equation} \label{nclemmaeq2}
\tau(\delta(a)b) = - \tau (\sigma (a) \delta (b)), \,\,\,\,\, \forall a,b \in A.
\end{equation}
To prove \eqref{eq:nclemmaeq}, first we assume that $m=1$. Now if $i=i_1=0$, then \eqref{eq:nclemmaeq} says that $\tau (ab) = \tau (b \sigma^{-1}(a))$ which is true because of the $\sigma$-trace property of $\tau$. By induction on $i$, assume that
\begin{equation} \label{nclemmaeqm=1}
\tau (b \sigma^{-1}(\delta \sigma^{-1})^{i}(a)) = (-1)^i \tau(a \delta^i(b)) \,\,\,\, \forall a,b \in A,
\end{equation}
for a fixed $i \geq 0$. Then by using \eqref{nclemmaeq0} and \eqref{nclemmaeq2}, we have
\begin{eqnarray}
\tau (b \sigma^{-1}(\delta \sigma^{-1})^{i+1}(a)) &=& \tau (b \sigma^{-1} (\delta \sigma^{-1})^{i} \delta \sigma^{-1}(a)) \nonumber \\
& =& (-1)^i \tau(\delta \sigma^{-1} (a) \delta^i(b)) \nonumber \\
&=& (-1)^{i+1} \tau (a \delta^{i+1}(b)). \nonumber
\end{eqnarray}
Therefore \eqref{nclemmaeqm=1} holds for all $i \geq 0$. Now to prove the general case, we use \eqref{nclemmaeqm=1} several times:
\begin{eqnarray}\
&& \tau \big ( b \sigma^{-1} (\delta \sigma^{-1})^{i_1} \cdots \sigma^{-1} (\delta \sigma^{-1})^{i_m} (a) \big ) \nonumber \\
&=& (-1)^{i_1} \tau \big ( \sigma^{-1} (\delta \sigma^{-1})^{i_2} \cdots \sigma^{-1} (\delta \sigma^{-1})^{i_m} (a) \delta^{i_1}(b) \big ) \nonumber \\
&=& (-1)^{i_1} \tau \big ( \sigma \delta^{i_1}(b) \sigma^{-1} (\delta \sigma^{-1})^{i_2} \cdots \sigma^{-1} (\delta \sigma^{-1})^{i_m} (a)  \big ) \nonumber \\
&=& (-1)^{i_1 + i_2} \tau \big ( \sigma^{-1} (\delta \sigma^{-1})^{i_3} \cdots \sigma^{-1} (\delta \sigma^{-1})^{i_m} (a) \delta^{i_2} \sigma \delta^{i_1}(b) \big ) \nonumber \\
&& \dots \nonumber \\
&=& (-1)^{i_1 + i_2 + \cdots + i_m} \tau \big( a \delta^{i_m} \sigma \delta^{i_{m-1}} \cdots \sigma \delta^{i_1} (b) \big). \nonumber
\end{eqnarray}

\end{proof}
\end{nctracelemma}

\newtheorem{nctrace}[twisteddef]{Theorem}
\begin{nctrace} \label{nctrace}

Let $A, \sigma, \delta, \tau$ be as in Definition \ref{twisteddef} where $\tau$ is a $\sigma$-trace and $\delta$ is $\sigma$-derivation. If $\tau \circ \delta = 0$, then the linear functional $\textnormal{res} \, : \Psi (A, \sigma , \delta)  \to \mathbb{C}$ defined by
\[ \textnormal{res} \, \, \big ( \sum_{i=- \infty}^{n} a_i \xi^i \big ) = \tau (a_{-1})  \]
is a trace.

\begin{proof}
Let $a,b \in A$, $m,n \in \mathbb{Z}$. We shall show that
\begin{equation} \label{identity}
\textnormal{res} \, (a \xi^n b \xi^m)=\textnormal{res} \, ( b \xi^m a \xi^n).
\end{equation}
One can easily see that both sides of \eqref{identity} are $0$ if $n,m \geq 0$, or if $n,m < 0$. So it suffices to prove the identity for $n \geq 0$, $m < 0$. In this case, by \eqref{eq:twistedrelation2} we have
\[a \xi^n b \xi^m = \sum_{i=0}^{n} a P_{i,\, n} (\sigma , \delta) (b) \xi^{i+m}.\]
Therefore by definition of $\textnormal{res} \,$, if $n+m < -1$, then $ \textnormal{res} \, (a \xi^n b \xi^m)=0$; otherwise we have:
\begin{equation} \label{eq:nctraceL}
 \textnormal{res} \, (a \xi^n b \xi^m) = \tau ( a P_{-m-1, \, n}(\sigma, \delta)(b)).
\end{equation}
Also by \eqref{eq:twistedrelation-n}
\[ b \xi^m a \xi^n = \sum_{i=(i_1, \dots , i_{-m}) \geq 0} (-1)^{| i |} b \sigma^{-1} (\delta \sigma^{-1})^{i_1} \cdots \sigma^{-1} (\delta \sigma^{-1})^{i_{-m}} (a) \xi^{m+n- |i|}. \]
Again by definition, if $n+m <-1$, $\textnormal{res} \, ( b \xi^m a \xi^n)=0$; otherwise by using Lemma \ref{nctrace}, we have:
\begin{eqnarray}
\textnormal{res} \, ( b \xi^m a \xi^n) &=& (-1)^{m+n+1} \tau \Big ( \sum_{|i|= m+n+1} b \sigma^{-1} (\delta \sigma^{-1})^{i_1} \cdots \sigma^{-1} (\delta \sigma^{-1})^{i_{-m}} (a) \Big ) \nonumber \\
&=&  \tau \Big( \sum_{|i|=m+n+1} a \delta^{i_{-m}} \sigma \delta^{i_{m-1}} \cdots \sigma \delta^{i_1} (b) \Big) \nonumber \\
&=& \tau \big( a P_{-m-1,\, n}(\sigma, \delta) (b) \big). \nonumber
\end{eqnarray}
Hence by \eqref{eq:nctraceL}, the identity \eqref{identity} is proved.

\end{proof}

\end{nctrace}

We notice that if $\delta \circ \sigma = \sigma \circ \delta$, the multiplication formulas \eqref{eq:twistedrelation2} and \eqref{eq:twistedrelation-n} in $\Psi(A,\sigma,\delta)$ simplify quite a bit, and reduce to the following:

\begin{equation} \label{commutativecommutation}
\xi^n a = \sum_{j=0}^\infty \binom{n}{j} \delta^j (\sigma^{n-j}(a))  \xi^{n-j}, \qquad n \in \mathbb{Z}.
\end{equation}

A special case of \eqref{commutativecommutation} is the algebra of $q$-pseudodifferential symbols on the circle defined in \cite{klr}. In this case $A=C^{\infty}(S^1), \sigma (f) = f(qx)$, and $\delta(f)= \frac{f(qx) - f(x) }{q-1}$. It is easy to see that $\delta \circ \sigma = \sigma \circ \delta$, and the resulting algebra coincides with $\Psi(A, \sigma, \delta).$

\newtheorem{traceproperty}[twisteddef]{Lemma}
\begin{traceproperty}\label{traceproperty}
Let $\sigma$ be an automorphism of an algebra $A$ and $\delta$ be a $\sigma$-derivation of $A$ such that $\delta \circ \sigma = \sigma \circ \delta.$

\begin{enumerate}
\item For any non-negative integer $k$, and $a,b \in A$, we have
\[\delta^k (ab)= \sum_{i=0}^{k} \binom{k}{i} \sigma^i \delta^{k-i}(a) \delta^i(b).\]

\item If $\varphi :A \to \mathbb{C}$ is a linear functional such that $\varphi \circ \delta = 0$, then for any non-negative integer $k$, and $a,b \in A$, we have
\[ \varphi (\delta^k(a)b) = (-1)^k \varphi (\sigma^k(a) \delta^k(b)). \]
\end{enumerate}

\begin{proof}
One can prove 1 by induction. To prove 2, since $\varphi(\delta(ab)) = \varphi (\delta(a) b + \sigma (a) \delta(b))=0$, we have
\[ \varphi(\delta(a)b) = - \varphi (\sigma (a) \delta (b)), \qquad \forall a,b \in A. \]
Now for any non-negative integer $k$, we have
\begin{eqnarray}
\varphi (\delta^k(a)b) &=& - \varphi (\sigma(\delta^{k-1}(a))\delta(b))= - \varphi (\delta^{k-1}(\sigma(a))\delta(b)) \nonumber \\
&=& \cdots = (-1)^k \varphi (\sigma^k(a) \delta^k(b)). \nonumber
\end{eqnarray}
\end{proof}

\end{traceproperty}

\newtheorem{twistedtrace}[twisteddef]{Proposition}

\begin{twistedtrace}\label{twistedtrace}
Let $A$ be an algebra, $\sigma$ an automorphism of $A$, and $\delta : A \to A$ a $\sigma$-derivation such that $\delta \circ \sigma = \sigma \circ \delta$.
\begin{enumerate}
\item The map $\sigma' : \Psi (A, \sigma, \delta) \to \Psi (A, \sigma, \delta)$ defined by
\[ \sigma'(\sum_{n=- \infty}^N a_n  \xi^{n}) = \sum_{n=- \infty}^N \sigma (a_n)  \xi^{n}, \]
is an automorphism of the algebra of formal twisted pseudodifferential symbols.
\item If $\tau$ is a trace on $A$ such that $\tau \circ \sigma =\tau $, and $\tau \circ \delta = 0$, then the map $ \textnormal{res}  : \Psi (A, \sigma, \delta) \to \mathbb{C} $ defined by
\[ \textnormal{res} \, (\sum_{n=- \infty}^N a_n  \xi^{n})= \tau (a_{-1})\]
is a $\sigma'^{-1}$-trace on $\Psi(A,\sigma,\delta)$.
\end{enumerate}

\begin{proof}

\begin{enumerate}
\item Obviously the map $\sigma'$ is a linear isomorphism. Also for any $a,b \in A$, and $m,n \in \mathbb{Z}$, we have
\begin{eqnarray}
\sigma' (a  \xi^n b  \xi^m) &=& \sigma' \Big ( \sum_{j=0}^\infty \binom{n}{j} a \delta^j (\sigma^{n-j}(b))  \xi^{n+m-j}\Big ) \nonumber \\
&=& \sum_{j=0}^\infty \binom{n}{j} \sigma (a) \sigma (\delta^j (\sigma^{n-j}(b)))  \xi^{n+m-j} \nonumber \\
&=& \sum_{j=0}^\infty \binom{n}{j} \sigma (a) \delta^j (\sigma^{n+1-j}(b))  \xi^{n+m-j}, \nonumber
\end{eqnarray}
and
\begin{eqnarray}
\sigma' (a  \xi^n ) \sigma'(b  \xi^m) &=& \sigma (a)  \xi^n  \sigma (b)  \xi^m  \nonumber \\
&=& \sum_{j=0}^\infty \binom{n}{j} \sigma (a) \delta^j (\sigma^{n-j+1}(b))  \xi^{n+m-j}. \nonumber
\end{eqnarray}
Therefore $\sigma'$ is an automorphism.

\item It suffices to show that for any $a,b \in A$, and $m,n \in \mathbb{Z}$,
\begin{equation} \label{twistedtraceeq}
\textnormal{res} \, (a  \xi^n \sigma (b)  \xi^m) = \textnormal{res} \, (b  \xi^m a  \xi^n  ).
\end{equation}

We have

\[ a  \xi^n \sigma (b)  \xi^m = \sum_{j=0}^\infty \binom{n}{j} a \delta^j (\sigma^{n-j+1}(b))  \xi^{n+m-j}\]
and
\[ b  \xi^m a  \xi^n = \sum_{j=0}^\infty \binom{m}{j} b \delta^j (\sigma^{m-j}(a))  \xi^{m+n-j}. \]
If $n+m < -1$, then  both sides of \eqref{twistedtraceeq} are $0$. If $n+m = -1$, then
\[ \textnormal{res} \, (a  \xi^n \sigma (b)  \xi^m) = \tau (a \sigma^{n+1}(b)),  \]
and
\begin{eqnarray}
\textnormal{res} \, (b  \xi^m a  \xi^n) &=& \tau (b \sigma^m (a)) = \tau (b \sigma^{-n-1} (a)) \nonumber \\
&=& \tau (\sigma^{n+1} (b) a) = \tau (a \sigma^{n+1} (b) ). \nonumber
\end{eqnarray}

Now assume that $n+m=k > -1$, then

\[ \textnormal{res} \, (a  \xi^n \sigma (b)  \xi^m) = \binom{n}{k+1} \tau (a \delta^{k+1}(\sigma^{n-k}(b))) \]
and
\begin{eqnarray}
\textnormal{res} \, (b  \xi^m a  \xi^n) &=&  \binom{m}{k+1}\tau (b \delta^{k+1}(\sigma^{m-k-1}(a))) \nonumber \\
&=& \binom{k-n}{k+1} \tau (b \delta^{k+1}(\sigma^{-n-1}(a))) \nonumber \\
&=& (-1)^{k+1} \binom{n}{k+1} \tau (b \delta^{k+1}(\sigma^{-n-1}(a))). \nonumber
\end{eqnarray}
Now by using Lemma \ref{traceproperty}, the desired result follows since
\begin{eqnarray}
\tau(a \delta^{k+1}(\sigma^{n-k}(b))) &=& \tau( \delta^{k+1}(\sigma^{n-k}(b)) a)  \nonumber \\
&=& (-1)^{k+1} \tau (\sigma^{n+1}(b) \delta^{k+1}(a)) \nonumber \\
&=& (-1)^{k+1} \tau (b \delta^{k+1}(\sigma^{-n-1}(a))). \nonumber
\end{eqnarray}

\end{enumerate}
\end{proof}
\end{twistedtrace}

\section{A multidimensional case}

In this section we consider an algebra $A$,  an automorphism $\sigma: A \to A $,
and two $\sigma$-derivations $\delta_1 , \delta_2 : A \to A$.
Let  $\Psi_0 (A, \sigma, \delta_1, \delta_2)$ be the  set of all
formal power series in noncommuting variables $\xi_1$ and $\xi_2$,
with coefficients in $A$, of the form
\[ D= \sum a_{i_1,j_1,\dots,i_m,j_m} \, \xi_1^{i_1} \xi_2^{j_1} \cdots \xi_1^{i_m} \xi_2^{j_m} \]
where $m= m(D)$ is a positive integer,
 and the summation is over all $m$-tuples of integers
 $i=(i_1, \dots, i_m)$ and $j=(j_1,\dots ,j_m)$ in $\mathbb{Z}^m$  such that $i_k$ and $j_k$ are less than some
 $N =N(D)\in \mathbb{Z}$,  for all $k=1, \dots, m$; and
 $a_{i_1,j_1,\dots,i_m,j_m} \in A$.
To define the multiplication, we impose the relations
\[
\xi_j a - \sigma (a) \xi_j = \delta_j (a),
\]
\[
\xi_j^{-1} \xi_j = \xi_j \xi_j^{-1}=1.
\]
for all $a \in A, j=1,2$.
Therefore by using the identities of Section 3, for any $a \in A$, integer $n \geq 1$, and $j=1,2,$ we have

\begin{equation} \label{eq:multitwistedrelation2}
\xi_j^n a = \sum_{i=0}^{n} P_{i, \, n}(\sigma , \delta_j)(a) \xi_j^{i}, \nonumber
\end{equation}
and
\begin{equation} \label{eq:multitwistedrelation-n}
\xi_j^{-n}a = \sum_{i_1=0}^{\infty} \cdots \sum_{i_{n}=0}^{\infty} (-1)^{i_1  + \cdots + i_{n}} \sigma^{-1} (\delta_j \sigma^{-1})^{i_{n}} \cdots \sigma^{-1} (\delta_j \sigma^{-1})^{i_1} (a) \xi_j^{-n - i_1 - \cdots - i_{n}}. \nonumber
\end{equation}

Now we define the algebra $\Psi (A, \sigma, \delta_1, \delta_2)$ to be the quotient of $\Psi_0 (A, \sigma, \delta_1, \delta_2)$ by the two sided ideal generated by $\xi_1 \xi_2 - \xi_2 \xi_1$. Note that each element of the latter has a representation of the form
\[ \sum_{i=-\infty}^{M} \sum_{j=-\infty}^{N} a_{ij} \, \xi_1^i \xi_2^j, \,\,\, a_{ij} \in A, \]
which is not necessarily unique.

\newtheorem{multitracelemma}[twisteddef]{Lemma}
\begin{multitracelemma}\label{multitracelemma}
Let $A, \sigma, \delta_1, \delta_2$ be as above, and $\tau : A \to \mathbb{C}$ be a $\sigma^2$-trace such that $\tau \circ \delta_j =0$ for $j=1,2$, and $\tau \circ \sigma = \tau .$ Then for any $a,b \in A$, and tuples of non-negative integers $i=(i_1, \dots , i_q )$, and $j=(j_1, \dots , j_p)$, we have
\begin{eqnarray}\label{multilemma1}
&& \tau \big ( b \sigma^{-1} (\delta_1 \sigma^{-1})^{j_1} \cdots \sigma^{-1} (\delta_1 \sigma^{-1})^{j_{p}} \sigma^{-1} (\delta_2 \sigma^{-1})^{i_1} \cdots \sigma^{-1} (\delta_2 \sigma^{-1})^{i_q}(a) \big ) \nonumber \\
&=& (-1)^{|i|+|j|} \tau \big ( a \delta_2^{i_q} \sigma \delta_2^{i_{q-1}} \cdots \sigma \delta_2^{i_1} \delta_1^{j_p} \sigma \delta_1^{j_{p-1}} \cdots \sigma \delta_1^{j_1}(b) \big );
\end{eqnarray}
and
\begin{eqnarray}\label{multilemma2}
&& \tau \big ( a \delta_1^{j_1} \sigma \delta_1^{j_2} \cdots \sigma \delta_1^{j_p} \sigma^{-1}(\delta_2 \sigma^{-1})^{i_1} \cdots \sigma^{-1}(\delta_2 \sigma^{-1})^{i_q}(b) \big ) \nonumber \\
&=& (-1)^{|i|+|j|} \tau \big ( b \delta_2^{i_q} \sigma \delta_2^{i_{q-1}} \cdots \sigma \delta_2^{i_1} \sigma^{-1} (\delta_1 \sigma^{-1})^{j_p} \cdots \sigma^{-1} (\delta_1 \sigma^{-1})^{j_1}(a) \big ).
\end{eqnarray}

\begin{proof}
First we show by induction that for $a,b \in A$, non-negative integer $i$, and $\delta = \delta_j, j=1,2$, we have
\begin{equation} \label{formultilemma1}
\tau \big ( b \sigma^{-1} (\delta \sigma^{-1})^i(a) \big ) = (-1)^i \tau \big ( \sigma(a) \delta^i(b) \big ).
\end{equation}
If $i=0$, then $\tau(b \sigma^{-1}(a)) = \tau (\sigma (a)b)$, because $\tau$ is a $\sigma^2$-trace. Now assume that \eqref{formultilemma1} holds for $i$. Then we have

\begin{eqnarray}
\tau \big ( b \sigma^{-1} (\delta \sigma^{-1})^{i+1}(a) \big ) &=& \tau \big ( b \sigma^{-1} (\delta \sigma^{-1})^i \delta \sigma^{-1} (a) \big ) = (-1)^i \tau \big ( \sigma \delta \sigma^{-1}(a) \delta^i(b) \big ) \nonumber \\
&=& (-1)^i \tau \big (  \delta \sigma^{-1}(a) \sigma^{-1} \delta^i(b) \big ) = (-1)^i \tau \big ( \sigma \delta^i(b) \delta \sigma^{-1}(a)  \big ) \nonumber \\
&=& (-1)^{i+1} \tau \big ( \delta^{i+1} (b) \sigma^{-1}(a)  \big ) = (-1)^{i+1} \tau \big ( \sigma (a) \delta^{i+1} (b)   \big ). \nonumber
\end{eqnarray}

Also one can see by induction that
\begin{equation} \label{formultilemma2}
\tau \big (  (\delta \sigma^{-1})^{i}(a) \, b  \big ) = (-1)^i \tau \big ( a \, \delta^i (b) \big );
\end{equation}
and
\begin{equation} \label{formultilemma3}
\tau \big ( a \, \delta^i (b) \big ) = (-1)^i \tau \big ( \sigma^2 (b) \, (\delta \sigma^{-1})^{i}(a) \big ).
\end{equation}

Now by using \eqref{formultilemma1} and \eqref{formultilemma2} several times, we can prove \eqref{multilemma1}:
\begin{eqnarray}
&& \tau \big ( b \sigma^{-1} (\delta_1 \sigma^{-1})^{j_1} \cdots \sigma^{-1} (\delta_1 \sigma^{-1})^{j_{p}} \sigma^{-1} (\delta_2 \sigma^{-1})^{i_1} \cdots \sigma^{-1} (\delta_2 \sigma^{-1})^{i_q}(a) \big ) \nonumber \\
&=& (-1)^{j_1} \tau \big ( (\delta_1 \sigma^{-1})^{j_2} \cdots \sigma^{-1} (\delta_1 \sigma^{-1})^{j_{p}}  \nonumber \\
&& \qquad \qquad \qquad \qquad \qquad \qquad  \sigma^{-1} (\delta_2 \sigma^{-1})^{i_1} \cdots \sigma^{-1} (\delta_2 \sigma^{-1})^{i_q}(a) \delta_1^{j_1}(b) \big ) \nonumber \\
&=& (-1)^{j_1} \tau \big (\sigma \delta_1^{j_1}(b) \sigma^{-1} (\delta_1 \sigma^{-1})^{j_2} \cdots \sigma^{-1} (\delta_1 \sigma^{-1})^{j_{p}} \nonumber \\
&& \qquad \qquad \qquad \qquad \qquad \qquad \qquad \sigma^{-1} (\delta_2 \sigma^{-1})^{i_1} \cdots \sigma^{-1} (\delta_2 \sigma^{-1})^{i_q}(a)  \big ) \nonumber \\
&=& \nonumber \\
&& \cdots \nonumber \\
&=& (-1)^{|j|} \tau \big (  \sigma \delta_1^{j_p} \sigma \delta_1^{j_{p-1}} \cdots \sigma \delta_1^{j_1}(b) \nonumber \\
&& \qquad \qquad \qquad \qquad \sigma^{-1} (\delta_2 \sigma^{-1})^{i_1} \sigma^{-1} (\delta_2 \sigma^{-1})^{i_2} \cdots \sigma^{-1} (\delta_2 \sigma^{-1})^{i_q} (a) \big ) \nonumber \\
&=& (-1)^{|j|} \tau \big ( (\delta_2 \sigma^{-1})^{i_1} \sigma^{-1} (\delta_2 \sigma^{-1})^{i_2} \cdots \sigma^{-1} (\delta_2 \sigma^{-1})^{i_q} (a) \delta_1^{j_p} \sigma \delta_1^{j_{p-1}} \cdots \sigma \delta_1^{j_1}(b) \big ) \nonumber \\
&=& (-1)^{|j|+i_1} \tau \big ( \sigma^{-1} (\delta_2 \sigma^{-1})^{i_2} \cdots \sigma^{-1} (\delta_2 \sigma^{-1})^{i_q} (a) \delta_2^{i_1} \delta_1^{j_p} \sigma \delta_1^{j_{p-1}} \cdots \sigma \delta_1^{j_1}(b) \big ) \nonumber \\
&=& (-1)^{|j|+i_1} \tau \big (  (\delta_2 \sigma^{-1})^{i_2} \cdots \sigma^{-1} (\delta_2 \sigma^{-1})^{i_q} (a) \sigma \delta_2^{i_1} \delta_1^{j_p} \sigma \delta_1^{j_{p-1}} \cdots \sigma \delta_1^{j_1}(b) \big ) \nonumber \\
&=& \nonumber \\
&& \cdots \nonumber \\
&=& (-1)^{|i|+|j|} \tau \big ( a \delta_2^{i_q} \sigma \delta_2^{i_{q-1}} \cdots \sigma \delta_2^{i_1} \delta_1^{j_p} \sigma \delta_1^{j_{p-1}} \cdots \sigma \delta_1^{j_1}(b) \big ). \nonumber
\end{eqnarray}

By using \eqref{formultilemma2} and \eqref{formultilemma3}, we prove \eqref{multilemma2}:
\begin{eqnarray}
&& \tau \big ( a \delta_1^{j_1} \sigma \delta_1^{j_2} \cdots \sigma \delta_1^{j_p} \sigma^{-1}(\delta_2 \sigma^{-1})^{i_1} \cdots \sigma^{-1}(\delta_2 \sigma^{-1})^{i_q}(b) \big ) \nonumber \\
&=& (-1)^{j_1} \tau \big (  \sigma^3 \delta_1^{j_2} \cdots \sigma \delta_1^{j_p} \sigma^{-1}(\delta_2 \sigma^{-1})^{i_1} \cdots \sigma^{-1}(\delta_2 \sigma^{-1})^{i_q}(b)  (\delta_1 \sigma^{-1})^{j_1} (a)      \big ) \nonumber \nonumber \\
&=& (-1)^{j_1} \tau \big ( \sigma^{-1} (\delta_1 \sigma^{-1})^{j_1} (a) \delta_1^{j_2} \cdots \sigma \delta_1^{j_p} \sigma^{-1}(\delta_2 \sigma^{-1})^{i_1} \cdots \sigma^{-1}(\delta_2 \sigma^{-1})^{i_q}(b)  \big ) \nonumber \\
&=& \nonumber \\
&& \cdots \nonumber \\
&=& (-1)^{|j|} \tau \big ( \sigma (\delta_2 \sigma^{-1})^{i_1} \sigma^{-1} (\delta_2 \sigma^{-1})^{i_2} \cdots \sigma^{-1} (\delta_2 \sigma^{-1})^{i_q}(b) \nonumber \\
&& \qquad \qquad \qquad \qquad \qquad (\delta_1 \sigma^{-1})^{j_p} \sigma^{-1} (\delta_1 \sigma^{-1})^{j_{p-1}} \cdots \sigma^{-1} (\delta_1 \sigma^{-1})^{j_1}(a) \big ) \nonumber
\end{eqnarray}
\begin{eqnarray}
&=& (-1)^{|j|} \tau \big ( (\delta_2 \sigma^{-1})^{i_1} \sigma^{-1} (\delta_2 \sigma^{-1})^{i_2} \cdots \sigma^{-1} (\delta_2 \sigma^{-1})^{i_q}(b) \nonumber \\
&& \qquad \qquad \qquad \qquad \qquad \sigma^{-1} (\delta_1 \sigma^{-1})^{j_p} \sigma^{-1} (\delta_1 \sigma^{-1})^{j_{p-1}} \cdots \sigma^{-1} (\delta_1 \sigma^{-1})^{j_1}(a) \big ) \nonumber \\
&=& (-1)^{|j|+i_1} \tau \big ( \sigma^{-1} (\delta_2 \sigma^{-1})^{i_2} \cdots \sigma^{-1} (\delta_2 \sigma^{-1})^{i_q}(b) \nonumber \\
&& \qquad \qquad \qquad \qquad \delta_2^{i_1} \sigma^{-1} (\delta_1 \sigma^{-1})^{j_p} \sigma^{-1} (\delta_1 \sigma^{-1})^{j_{p-1}} \cdots \sigma^{-1} (\delta_1 \sigma^{-1})^{j_1}(a) \big ) \nonumber \\
&=& (-1)^{|j|+i_1} \tau \big (  (\delta_2 \sigma^{-1})^{i_2} \cdots \sigma^{-1} (\delta_2 \sigma^{-1})^{i_q}(b) \nonumber \\
&& \qquad \qquad \qquad \qquad \sigma \delta_2^{i_1} \sigma^{-1} (\delta_1 \sigma^{-1})^{j_p} \sigma^{-1} (\delta_1 \sigma^{-1})^{j_{p-1}} \cdots \sigma^{-1} (\delta_1 \sigma^{-1})^{j_1}(a) \big ) \nonumber \\
&=& \nonumber \\
&& \cdots \nonumber \\
&=& (-1)^{|i|+|j|} \tau \big ( b \delta_2^{i_q} \sigma \delta_2^{i_{q-1}} \cdots \sigma \delta_2^{i_1} \sigma^{-1} (\delta_1 \sigma^{-1})^{j_p} \cdots \sigma^{-1} (\delta_1 \sigma^{-1})^{j_1}(a) \big ). \nonumber
\end{eqnarray}

\end{proof}

\end{multitracelemma}

Now we can define the $2$-dimensional analogue of the twisted residue map.

\newtheorem{multitracethm}[twisteddef]{Theorem}
\begin{multitracethm}\label{multitracethm}
Let $A, \sigma, \delta_1, \delta_2, \tau$ be as in Lemma \ref{multitracelemma}, where $\tau$ is an invariant $\sigma^2$-trace on $A$; and define a linear map $\textnormal{res}_0 \, : \Psi_0 (A, \sigma, \delta_1, \delta_2) \to \mathbb{C}$ by
\[ \textnormal{res} _0  \, \big ( \sum a_{i_1,j_1,\dots,i_m,j_m} \, \xi_1^{i_1} \xi_2^{j_1} \cdots \xi_1^{i_m} \xi_2^{j_m} \big ) =
\sum_{|i| = -1, |j|= -1} \tau (a_{i_1,j_1,\dots,i_m,j_m}), \]
where $i=(i_1, \dots , i_m)$ and $j=(j_1, \dots , j_m).$ If $ (\xi_1 \xi_2 - \xi_2 \xi_1) \subset \textnormal{Ker} (\textnormal{res}_0)$, then the induced map $\textnormal{res} : \Psi (A, \sigma, \delta_1, \delta_2) \to \mathbb{C}$ given by
\[ \textnormal{res} \, \big ( \sum_{i=-\infty}^{M} \sum_{j=-\infty}^{N} a_{ij} \, \xi_1^i \xi_2^j    \big ) = \tau (a_{-1 -1}),\]
is a trace functional.

\begin{proof}
Let $a,b \in A$, and $m,n,p,q \in \mathbb{Z}.$ We shall show that
\begin{equation} \label{eq:multitraceeq}
\textnormal{res} \, (a \, \xi_1^m \xi_2^n \, b \, \xi_1^p \xi_2^q) = \textnormal{res} \, (b \, \xi_1^p \xi_2^q \, a \, \xi_1^m \xi_2^n ).
\end{equation}
One can easily see that if $m,p \geq 0$, or if $n,q \geq 0$, then by definition of $\textnormal{res} \,$, both sides of \eqref{eq:multitraceeq} are 0. So it suffices to prove the identity for two general cases:
\begin{enumerate}
\item $m,n \geq 0$, and $p,q <0.$ In this case
\begin{eqnarray}
a \, \xi_1^m \xi_2^n \, b \, \xi_1^p \xi_2^q &=& a \, \xi_2^n \xi_1^m  \, b \, \xi_1^p \xi_2^q \nonumber \\
&=& \sum_{j=0}^{m} \sum_{i=0}^{n} a P_{i,n}(\sigma, \delta_2) P_{j,m}(\sigma, \delta_1)(b) \xi_1^{p+j} \xi_2^{q+i}. \nonumber
\end{eqnarray}
Therefore by definition of $\textnormal{res} \,$, if $n+q < -1$ or $m+p < -1$, then $\textnormal{res} \, (a \, \xi_1^m \xi_2^n \, b \, \xi_1^p \xi_2^q)=0$, otherwise we have
\begin{eqnarray} \label{multitracel}
&& \textnormal{res} \, (a \, \xi_1^m \xi_2^n \, b \, \xi_1^p \xi_2^q) = \tau \big ( a P_{-q-1,n}(\sigma, \delta_2) P_{-p-1,m} (\sigma, \delta_1) (b) \big ).
\end{eqnarray}
On the other hand we have
\begin{eqnarray}
&& b \, \xi_1^p \xi_2^q \, a \, \xi_1^m \xi_2^n \nonumber \\
&=& \sum (-1)^{|i|+|j|} b \sigma^{-1} (\delta_1 \sigma^{-1})^{j_1} \cdots \sigma^{-1} (\delta_1 \sigma^{-1})^{j_{-p}} \nonumber \\
&& \qquad \qquad \sigma^{-1} (\delta_2 \sigma^{-1})^{i_1} \cdots \sigma^{-1} (\delta_2 \sigma^{-1})^{i_{-q}} (a) \xi_1^{m+p-|j|} \xi_2^{n+q-|i|}, \nonumber
\end{eqnarray}
where the sum is over all tuples of non-negative integers $i=(i_1, \dots , i_{-q})$ and $j=(j_1, \dots , j_{-p})$.
Therefore by definition of $\textnormal{res} \,$, if $n+q < -1$ or $m+p < -1$, then $\textnormal{res} \, (b \, \xi_1^p \xi_2^q \, a \, \xi_1^m \xi_2^n) = 0$, otherwise by Lemma \ref{multitracelemma}
\begin{eqnarray} \label{multitracer}
&& \textnormal{res} \, (b \, \xi_1^p \xi_2^q \, a \, \xi_1^m \xi_2^n) \nonumber \\
&=& (-1)^{n+q+m+p} \sum \tau \big ( b \sigma^{-1} (\delta_1 \sigma^{-1})^{j_1} \cdots \sigma^{-1} (\delta_1 \sigma^{-1})^{j_{-p}}    \nonumber \\
&& \qquad \qquad \qquad \qquad \qquad \qquad \sigma^{-1} (\delta_2 \sigma^{-1})^{i_1} \cdots \sigma^{-1} (\delta_2 \sigma^{-1})^{i_{-q}} (a) \big )  \nonumber \\
&=& \sum \tau ( a \delta_2^{i_{-q}} \sigma \delta_2^{i_{-q-1}} \dots \sigma \delta_2^{i_1} \delta_1^{j_{-p}} \sigma \delta_1^{i_{-p-1}} \dots \sigma \delta_1^{j_1} (b) )
\end{eqnarray}
where the sums are over all tuples of non-negative integers $i=(i_1, \dots , i_{-q})$ and $j=(j_1, \dots , j_{-p})$, such that $|i| = n+q+1$, and $|j| = m+p+1.$ Therefore \eqref{multitracer} amounts to saying that
\begin{eqnarray}
&& \textnormal{res} \, (b \, \xi_1^p \xi_2^q \, a \, \xi_1^m \xi_2^n) = \tau \big ( a P_{-q-1,n}(\sigma, \delta_2) P_{-p-1, m}(\sigma, \delta_1)(b) \big ), \nonumber
\end{eqnarray}
which is equal to $\textnormal{res} \, (a \, \xi_1^m \xi_2^n \, b \, \xi_1^p \xi_2^q)$ by \eqref{multitracel}.
\item $m \geq 0, n < 0$, and $p<0, q \geq 0.$ In this case we have
\begin{eqnarray}
&& a \, \xi_1^m \xi_2^n \, b \, \xi_1^p \xi_2^q \nonumber \\
&=& \sum \sum_{j=0}^{m} (-1)^{|i|} a P_{j,m}(\sigma, \delta_1) \sigma^{-1} (\delta_2 \sigma^{-1})^{i_1} \nonumber \\
&& \qquad \qquad \qquad \qquad \cdots \sigma^{-1} (\delta_2 \sigma^{-1})^{i_{-n}}(b) \xi_1^{j+p} \xi_2^{n - |i|+q}, \nonumber
\end{eqnarray}
where the first summation is over all tuples of non-negative integers $i=(i_1, \dots , i_{-n}).$ Therefore by definition of $\textnormal{res} \,$, if $m+p < -1$, or $n+q < -1 $, then $\textnormal{res} \, (a \, \xi_1^m \xi_2^n \, b \, \xi_1^p \xi_2^q) = 0$, otherwise by Lemma \ref{multitracelemma} we have
\begin{eqnarray} \label{multitracel2}
&& \textnormal{res} \, (a \, \xi_1^m \xi_2^n \, b \, \xi_1^p \xi_2^q) \nonumber \\
&=& (-1)^{n+q+1} \sum_{|i|= n+q+1} \tau \big (a P_{-p-1, m}(\sigma, \delta_1)\sigma^{-1}(\delta_2 \sigma^{-1})^{i_1} \nonumber \\
&& \qquad \qquad \qquad \qquad \qquad \qquad \qquad \qquad \cdots \sigma^{-1}(\delta_2 \sigma^{-1})^{i_{-n}} (b) \big ) \nonumber \\
&=& (-1)^{n+q+1} \sum \tau \big ( a \delta_1^{r_1} \sigma \delta_1^{r_2} \dots \sigma \delta_1^{r_{-p}} \nonumber \\
&& \qquad \qquad \qquad \qquad \qquad  \sigma^{-1}(\delta_2 \sigma^{-1})^{i_1} \cdots \sigma^{-1}(\delta_2 \sigma^{-1})^{i_{-n}} (b) \big )  \nonumber \\
&=&(-1)^{m+p+1} \sum \tau \big ( b \delta_2^{i_{-n}} \sigma \delta_2^{i_{-n-1}} \cdots \sigma \delta_2^{i_1} \nonumber \\
&& \qquad \qquad \qquad \qquad \qquad \sigma^{-1} (\delta_1 \sigma^{-1})^{r_{-p}} \cdots \sigma^{-1} (\delta_1 \sigma^{-1})^{r_1}(a) \big ) \nonumber \\
&&
\end{eqnarray}
where in the last two sums, the summation is over all tuples of non-negative integers $i=(i_1, \dots , i_{-n})$ and $r=(r_1, \dots , r_{-p})$ such that $|i|=n+q+1$, and $|r|= m+p+1$.\\
On the other hand we have
\begin{eqnarray}
b \, \xi_1^p \xi_2^q \, a \, \xi_1^m \xi_2^n &=& b \, \xi_2^q \xi_1^p  \, a \, \xi_1^m \xi_2^n \nonumber \\
&=& \sum \sum_{k=0}^{q} (-1)^{|r|} b P_{k,\, q}(\sigma, \delta_2) \sigma^{-1} (\delta_1 \sigma^{-1})^{r_1} \nonumber \\
&& \qquad \qquad \qquad \qquad \cdots \sigma^{-1} (\delta_1 \sigma^{-1})^{r_{-p}}(a) \xi_1^{p-|r|+m} \xi_2^{n+k}, \nonumber
\end{eqnarray}
where the first summation is over all tuples of non-negative integers $r=(r_1, \dots , r_{-p})$. Therefore we have
\begin{eqnarray} \label{multitracer2}
&& \textnormal{res} \, (b \, \xi_1^p \xi_2^q \, a \, \xi_1^m \xi_2^n) \nonumber \\
&=& (-1)^{p+m+1} \sum \tau \big ( b P_{-n-1, \, q} (\sigma, \delta_2) \sigma^{-1} (\delta_1 \sigma^{-1})^{r_1} \nonumber \\
&& \qquad \qquad \qquad \qquad \qquad \qquad \qquad \cdots \sigma^{-1} (\delta_1 \sigma^{-1})^{r_{-p}}(a) \big ),
\end{eqnarray}
where the sum is over all tuples of non-negative integers $r=(r_1, \dots , r_{-p})$, such that $|r|=p+m+1$. \\

Hence $\textnormal{res} \,(a \, \xi_1^m \xi_2^n \, b \, \xi_1^p \xi_2^q) = \textnormal{res} \, (b \, \xi_1^p \xi_2^q \, a \, \xi_1^m \xi_2^n)$ by \eqref{multitracel2}, \eqref{multitracer2}, and the definition of $P_{-n-1, \, q}(\sigma, \delta_2)$.
\end{enumerate}

\end{proof}

\end{multitracethm}

In the following proposition, we show that if the twisted derivations commute with each other, then the condition $ (\xi_1 \xi_2 - \xi_2 \xi_1) \subset \textnormal{Ker} (\textnormal{res}_0)$ of Theorem \ref{multitracethm} holds.

\newtheorem{criterion}[twisteddef]{Proposition}
\begin{criterion}
Let $A, \sigma, \delta_1, \delta_2, \tau$ be as in Theorem \ref{multitracethm}, where $\tau$ is an invariant $\sigma^2$-trace on $A$. If $\delta_1 \circ \delta_2 = \delta_2 \circ \delta_1$, then $ (\xi_1 \xi_2 - \xi_2 \xi_1) \subset \textnormal{Ker} (\textnormal{res}_0)$.

\begin{proof}

 Let $a,b \in A$, and $j=1$ or $2$, we have $\tau (a \, \sigma \delta_j (b)) = \tau (a \, \delta_j \sigma (b))$, because:
\begin{eqnarray}
\tau (a \, \sigma \delta_j (b)) &=& \tau(\sigma^{-1}(a) \, \delta_j(b)) = - \tau(\delta_j \sigma^{-2}(a) \, b) \nonumber \\
&=& - \tau (\sigma^2 (b) \, \delta_j \sigma^{-2}(a)) = \tau (\delta_j \sigma (b) \, \sigma^{-2}(a)) \nonumber \\
&=&  \tau (a \, \delta_j \sigma (b)). \nonumber
\end{eqnarray}
Similarly one can see that $\tau (a \, \sigma^{-1} \delta_j (b)) = \tau (a \, \delta_j \sigma^{-1} (b))$. Therefore
\[ \tau (a\, \sigma^{n_1}\delta_1^{i_1}\delta_2^{j_2} \cdots \sigma^{n_r}\delta_1^{i_r}\delta_2^{j_r} (b)) \]
depends only on $|n|,|i|$ and $|j|$, for any $n=(n_1, \dots , n_r)$ an $r$-tuple of integers, and $i=(i_1,\dots, i_r)$ and $j=(j_1, \dots, j_r)$ any $r$-tuples of non-negative integers.\\

By using the above fact, one can see that for any choice of integers $i_k, j_k, i'_k, j'_k$,
\begin{eqnarray}
&& \textnormal{res}_0 \big ( a \, \xi_1^{i_1} \xi_2^{j_1} \cdots \xi_1^{i_m} \xi_2^{j_m} (\xi_1 \xi_2 - \xi_2 \xi_1) b \, \xi_1^{i'_1} \xi_2^{j'_1} \cdots \xi_1^{i'_{m'}} \xi_2^{j'_{m'}} \big ) \nonumber \\
&=& \textnormal{res}_0  \bigg ( a \, \xi_1^{i_1} \xi_2^{j_1} \cdots \xi_1^{i_m} \xi_2^{j_m} \Big ( \sigma^2 (b) (\xi_1 \xi_2 - \xi_2 \xi_1) + (\delta_1 \sigma(b) - \sigma \delta_1 (b)) \xi_2 + \nonumber \\
&& \qquad \qquad \qquad  (\sigma \delta_2 (b) - \delta_2 \sigma(b))\xi_1 + \delta_1 \delta_2 (b) - \delta_2 \delta_1 (b) \Big ) \xi_1^{i'_1} \xi_2^{j'_1} \cdots \xi_1^{i'_{m'}} \xi_2^{j'_{m'}} \bigg ) \nonumber \\
&=& 0. \nonumber 
\end{eqnarray}

\end{proof}

\end{criterion}

By a similar argument, one can show that starting with a $\sigma^n$-trace $\tau$ such that $\tau \circ \sigma = \tau$, and $n$ commuting $\sigma$-derivations $\delta_1, \dots , \delta_n : A \to A$ which leave $\tau$ invariant, the linear functional $\textnormal{res} : \Psi (A,\sigma, \delta_1, \dots , \delta_n) \to \mathbb{C}$ defined by $\textnormal{res} \, \Big ( \sum a_{i_1, \dots ,i_n} \, \xi_1^{i_1} \cdots \xi_n^{i_n} \Big )= \tau (a_{-1, \dots ,-1})$ is a trace. Here the algebra of twisted pseudodifferential symbols $\Psi (A,\sigma, \delta_1, \dots , \delta_n) $  is defined as in the $2$-dimensional case. We record this in:

\newtheorem{multitracethm-n}[twisteddef]{Theorem}
\begin{multitracethm-n}
Let $\sigma$ be an automorphism of an algebra $A$, and $\delta_1, \dots , \delta_n : A \to A$ commuting $\sigma$-derivations. If $\tau : A \to \mathbb{C}$ is a $\sigma^n$-trace such that $\tau \circ \sigma = \tau$, and $\tau \circ \delta_j = 0$ for $j=1, \dots , n$, then the linear functional
$\textnormal{res} : \Psi (A,\sigma, \delta_1, \dots , \delta_n) \to \mathbb{C}$ defined by
\[  \textnormal{res} \, \Big ( \sum a_{i_1, \dots ,i_n} \, \xi_1^{i_1} \cdots \xi_n^{i_n} \Big )= \tau (a_{-1, \dots ,-1})  \]
is a trace.
\end{multitracethm-n}

\section{The logarithmic cocycle}

In this section we extend the logarithmic cocycle of Kravchenko-Khesin \cite{krakhe} to
 our twisted setup. A special case of our construction is the logarithmic cocycle
 of Khesin-Lyubashenko-Roger \cite{klr} on the algebra of $q$-pseudodifferential
 symbols on the circle. In Section 3 we noticed that when the automorphism $\sigma : A \to A$ commutes
 with the $\sigma$-derivation $\delta$, the multiplication rules in $\Psi(A,\sigma,\delta)$ are
 derived from \eqref{commutativecommutation}. In this section we assume that $\delta$ and $\sigma$
 commute. To derive the main formulas \eqref{log}
 and \eqref{log1} of this section, in a formal manner we first assume that the algebra $A$ has
 a 1-parameter group of automorphisms $\{ \sigma_{t} \}_{t \in \mathbb{R}}$ such
 that $\sigma_1 = \sigma$. We shall then observe that Proposition \ref{logderivation}
  holds in general, without assuming the existence of $\sigma_t$. 

Let $\sigma_t$ be a 1-parameter group of automorphisms of $A$ with $\sigma_1 = \sigma$.
We can define an algebra of twisted pseudodifferential symbols with elements of
the form $\sum_{i=0}^{\infty} a_i \xi^{t-i}, \, t \in \mathbb{R}$. We can replace
the integer $n$ in \eqref{commutativecommutation}  by $t \in \mathbb{R}$ and obtain
\[
\xi^{t} a - \sigma_{t}(a) \xi^t = \sum_{j=1}^\infty \binom{t}{j} \delta^j (\sigma_{t-j}(a))  \xi^{t-j},
\]
from which the multiplication formula is derived.
By differentiating the above formula in $t$ at $t = 0$, and using the identity
\[ \frac{d}{dt} \mid_{t = 0} \xi^t = \log \xi \cdot \xi^t \mid_{t = 0} \, = \log \xi  , \]
we will have the commutation relation
\begin{equation} \label{log0}
[\log \xi, a] = \frac{d}{dt} \mid_{t = 0} \sigma_t (a) + \sum_{j=1}^{\infty} \frac{(-1)^{j-1}}{j} \sigma^{-j} \delta^j (a) \xi^{-j}.
\end{equation}
One would expect \eqref{log0} to define a derivation $[\log \xi \, , -] :  A \to \Psi(A,\sigma,\xi)$.
To define this  derivation, by assuming  $[\log \xi, \xi]=0$, we simply let
\begin{equation} \label{log}
[\log \xi, a \xi^n] = \frac{d}{dt} \mid_{t = 0} \sigma_t (a) \xi^n + \sum_{j=1}^{\infty} \frac{(-1)^{j-1}}{j} \sigma^{-j} \delta^j (a) \xi^{n-j},
\end{equation}
for any $a \in A, n\in \mathbb{Z}.$

\newtheorem{loglemma}[twisteddef]{Lemma}
\begin{loglemma}\label{loglemma}
Let $s \geq 1$ and $k \geq 0$ be integers and $\alpha \in \mathbb{R}$. Then

\begin{equation} \label{loglemmaid}
\frac{(-1)^{s-1}}{s} \binom{\alpha - s}{k} = \sum_{j=s}^{s+k} \frac{(-1)^{j-1}}{j} \binom{\alpha}{s+k-j} \binom{j}{j-s}.
\end{equation}

\begin{proof}
One can write \eqref{loglemmaid} in this form:
\[\binom{\alpha-s}{k}= \sum_{j=s}^{s+k} (-1)^{s+j} \binom{\alpha}{s+k-j} \binom{j-1}{s-1}.  \]
To prove the latter, we write:
\begin{eqnarray}
(1+x)^{\alpha-s} &=& (1+x)^\alpha(1+x)^{-s} \nonumber \\
&=& (1+x)^\alpha \Big ( 1 + \binom{-s}{1} x + \binom{-s}{2} x^2 + \binom{-s}{3} x^3 + \cdots \Big ) \nonumber \\
&=& (1+x)^\alpha + \binom{-s}{1} x (1+x)^\alpha + \binom{-s}{2} x^2 (1+x)^\alpha  + \cdots \nonumber \\
&=& (1+x)^\alpha - \binom{s}{1} x (1+x)^\alpha + \binom{s+1}{2} x^2 (1+x)^\alpha  + \cdots \nonumber
\end{eqnarray}
and consider the coefficient of $x^k$ on both sides.
\end{proof}

\end{loglemma}

For the following proposition,
 we assume that   $\sigma$ is an automorphism of $A$, and  $\sigma_t$ is
  a 1-parameter group of  automorphisms of $A$   such that
\[\frac{d}{dt} \mid_{t = 0} \sigma_t (\sigma(a)) = \sigma \big ( \frac{d}{dt} \mid_{t = 0} \sigma_t (a) \big ) \]
and
\[ \frac{d}{dt} \mid_{t = 0} \sigma_t (\delta(a)) = \delta \big (  \frac{d}{dt} \mid_{t = 0} \sigma_t (a) \big ) \]
for all $a \in A$.

\newtheorem{logderivation}[twisteddef]{Proposition}
\begin{logderivation}\label{logderivation}
The map $[\log \xi, -] : \Psi(A,\sigma,\delta) \to \Psi(A,\sigma,\delta)$ defined by \eqref{log} is a derivation.

\begin{proof}
Let $a,b \in A , n,m \in \mathbb{Z}.$ We have
\begin{eqnarray} \label{logderleft}
&& [\log \xi, a \xi^n b \xi^m] \nonumber \\
&=& [\log \xi, \sum_{k=0}^{\infty} \binom{n}{k} a \sigma^{n-k} \delta^k(b) \xi^{n-k+m}  ] \nonumber \\
&=& \sum_{k=0}^{\infty} \binom{n}{k} \frac{d}{dt} |_{t=0} \sigma_t (a \sigma^{n-k} \delta^k(b) ) \xi^{n+m-k} + \nonumber \\
&& \sum_{k=0}^{\infty} \sum_{j=1}^{\infty} \binom{n}{k} \frac{(-1)^{j-1}}{j} \sigma^{-j} \delta^j (a \sigma^{n-k} \delta^k(b)) \xi^{n+m-j-k}.
\end{eqnarray}

On the other hand we have:
\begin{eqnarray}
&& [\log \xi, a\xi^n] b\xi^m + a\xi^n [\log \xi, b \xi^m] \nonumber \\
&=&\Big ( \frac{d}{dt} \mid_{t = 0} \sigma_t (a) \xi^n + \sum_{j=1}^{\infty} \frac{(-1)^{j-1}}{j} \sigma^{-j} \delta^j (a) \xi^{n-j} \Big ) b \xi^m + \nonumber \\
&& a \xi^n \Big ( \frac{d}{dt} \mid_{t = 0} \sigma_t (b) \xi^m + \sum_{j=1}^{\infty} \frac{(-1)^{j-1}}{j} \sigma^{-j} \delta^j (b) \xi^{m-j} \Big ) \nonumber 
\end{eqnarray}
\begin{eqnarray}\label{logderright}
&=& \frac{d}{dt} \mid_{t = 0} \sigma_t (a) \sum_{k=0}^{\infty} \binom{n}{k} \sigma^{n-k} \delta^k (b) \xi^{n-k+m} + \nonumber \\
&& \sum_{j=1}^{\infty} \sum_{k=0}^{\infty}  \frac{(-1)^{j-1}}{j} \binom{n-j}{k} \sigma^{-j} \delta^j(a) \sigma^{n-j-k} \delta^k(b) \xi^{n-j-k+m} + \nonumber \\
&& \sum_{k=0}^{\infty} \binom{n}{k} a \sigma^{n-k} \delta^k (\frac{d}{dt} \mid_{t = 0} \sigma_t (b)) \xi^{n-k+m} + \nonumber \\
&& \sum_{j=1}^{\infty} \sum_{k=0}^{\infty} \frac{(-1)^{j-1}}{j} \binom{n}{k} a \sigma^{n-k-j} \delta^{k+j} (b) \xi^{n+m-k-j }.
\end{eqnarray}

Therefore the coefficient of $\xi^{n+m}$ in \eqref{logderleft} is
\begin{eqnarray}
&& \frac{d}{dt} \mid_{t = 0} \sigma _t (a \sigma^n (b)) \nonumber \\
&=& \frac{d}{dt} \mid_{t = 0} \sigma_t (a) \, \sigma^n(b) + a \frac{d}{dt} \mid_{t = 0} \sigma_t (\sigma^n (b))  , \nonumber
\end{eqnarray}
and the corresponding coefficient in \eqref{logderright} is
\begin{eqnarray}
\frac{d}{dt} \mid_{t = 0} \sigma_t (a) \, \sigma^n(b) + a \, \sigma^n \Big ( \frac{d}{dt} \mid_{t = 0} (b) \Big ). \nonumber
\end{eqnarray}

For any integer $r\geq 1$, the coefficient of $\xi^{n+m-r}$ in \eqref{logderleft} is
\[ \binom{n}{r} \frac{d}{dt} \mid_{t = 0} \sigma_t (a \sigma^{n-r} \delta^r (b)) + \sum \binom{n}{k} \frac{(-1)^{j-1}}{j} \sigma^{-j} \delta^j (a \sigma^{n-k} \delta^k(b))  \]
where the summation is over all integers $j\geq 1, k\geq 0$, such that $j+k =r$. By using Lemma \ref{traceproperty} and the derivation property of the derivative, the latter is equal to
\begin{eqnarray}\label{logderleft1}
&& \binom{n}{r} \frac{d}{dt} \mid_{t = 0} \sigma_t (a) \, \sigma^{n-r}\delta^r (b) + \binom{n}{r} a  \frac{d}{dt} \mid_{t = 0} \sigma_t ( \sigma^{n-r}\delta^r (b)) + \nonumber \\
&& \sum \sum_{i=0}^{j} \frac{(-1)^{j-1}}{j} \binom{n}{k} \binom{j}{i} \sigma^{i-j} \delta^{j-i}(a) \sigma^{n-r} \delta^{i+k}(b) \nonumber \\
&=& \binom{n}{r} \frac{d}{dt} \mid_{t = 0} \sigma_t (a) \, \sigma^{n-r}\delta^r (b) + \binom{n}{r} a  \frac{d}{dt} \mid_{t = 0} \sigma_t ( \sigma^{n-r}\delta^r (b)) + \nonumber \\
&& \sum \sum_{i=0}^{j-1} \frac{(-1)^{j-1}}{j} \binom{n}{k} \binom{j}{i} \sigma^{i-j} \delta^{j-i}(a) \sigma^{n-r} \delta^{i+k}(b) + \nonumber \\
&&  \sum \frac{(-1)^{j-1}}{j} \binom{n}{k} a \sigma^{n-r} \delta^{r}(b).
\end{eqnarray}
where both summations are over all integers $j\geq 1, k\geq 0$, such that $j+k =r$.
Also for any integer $r\geq 1$, the coefficient of $\xi^{n+m-r}$ in \eqref{logderright} is
\begin{eqnarray} \label{logderright1}
&& \binom{n}{r} \frac{d}{dt} \mid_{t = 0} \sigma_t (a) \, \sigma^{n-r} \delta^r (b) + \binom{n}{r} a \, \sigma^{n-r} \delta^r \Big ( \frac{d}{dt} \mid_{t = 0} \sigma_t (b) \Big ) + \nonumber \\
&&  \sum \frac{(-1)^{s-1}}{s} \binom{n-s}{l} \sigma^{-s} \delta^s (a) \sigma^{n-r} \delta^l(b) + \nonumber \\
&&  \sum \frac{(-1)^{s-1}}{s} \binom{n}{l} a \, \sigma^{n-r} \delta^{r}(b)
\end{eqnarray}
where both summations are over all integers $s\geq 1, l\geq 0$, such that $s+l =r$.

Therefore in order to show that \eqref{logderleft1} and \eqref{logderright1} are equal, it suffices to show that for fixed integers $s \geq 1, l \geq 0$, such that $s+l=r$:

\[ \frac{(-1)^{s-1}}{s} \binom{n-s}{l} = \sum \frac{(-1)^{j-1}}{j} \binom{n}{k} \binom{j}{i} \]
where the sum is over all integers $k \geq 0, j \geq 1, i \geq 0$ such that $i= j-s, i+k=l$, which amounts to saying that
\[\frac{(-1)^{s-1}}{s} \binom{n-s}{l} = \sum_{j=s}^{l+s} \frac{(-1)^{j-1}}{j} \binom{n}{l+s-j} \binom{j}{j-s}, \]
and the latter follows from Lemma \ref{loglemma}.

\end{proof}

\end{logderivation}

A close look at the proof of Proposition \ref{logderivation} shows that to define the derivation $[\log \xi, -]$, in fact we don't need to start with a one parameter group of automorphisms $\sigma_t$. All we need is a derivation $\delta_0 : A \to A$ such that $ \delta_0 \circ \sigma = \sigma \circ \delta_0$ and $\delta_0 \circ \delta = \delta \circ \delta_0$. We can then define

\begin{equation} \label{log1}
[\log \xi, a \xi^n] = \delta_0 (a) \, \xi^n + \sum_{j=1}^{\infty} \frac{(-1)^{j-1}}{j} \sigma^{-j} \delta^j (a) \xi^{n-j}.
\end{equation}
We find the following analogue of Proposition \ref{logderivation}.

\newtheorem{logderivation1}[twisteddef]{Proposition}
\begin{logderivation1}\label{logderivation1}
Let $\delta_0 : A \to A$ be a derivation such that $\delta_0 \circ \sigma = \sigma \circ \delta_0$ and $\delta_0 \circ \delta = \delta \circ \delta_0$. Then the map $[\log \xi, -] : \Psi(A,\sigma,\delta) \to \Psi(A,\sigma,\delta)$ defined by \eqref{log1} is a derivation.
\end{logderivation1}

Let $\tau : A \to \mathbb{C}$ be a $\delta$-invariant $\sigma$-trace on $A$. Using Theorem \ref{nctrace}, we have a trace $\textnormal{res} : \Psi(A, \sigma, \delta) \to \mathbb{C}.$ Using the derivation property of $[\log \xi, -]$ and the trace property of $\textnormal{res}$, it follows that the $2$-cochain
\[ c \, (D_1, D_2) = \textnormal{res} \, ([\log \xi, D_1] D_2)  \]
is a Lie algebra $2$-cocycle on $\Psi (A, \sigma, \delta).$ This extends the Radul $2$-cocycle \cite{krakhe} and its $q$-analogue in \cite{klr}.

\section{Twisting by cocycles}
In this section, using group 1-cocycles,  we give a general method to construct twisted
derivations, twisted traces, and twisted connections on a crossed
product algebra. Using twisted connections, one can then construct twisted spectral triples in
the sense of Connes-Moscovici \cite{conmos}. 
As an example, we recover  a twisted spectral
triple first constructed in \cite{conmos}.

Let $A$ be an algebra with a right action of a group $\Gamma$ by
automorphisms:
\[ A \times \Gamma \to A,  \qquad
 (a,\gamma) \mapsto a \cdot \gamma .\]
We consider the algebraic crossed product $A \rtimes \Gamma$ with the standard multiplication:
\[ (a \otimes \gamma)(b \otimes \mu) = (a \cdot \mu) b \otimes \gamma \mu, \quad \,\,\, a,b \in A, \,\,\, \gamma, \mu \in \Gamma. \]
Let $Z(A)$ denote the center of the algebra $A$, and $A^*$ its group of invertible elements.

\newtheorem{chain}[twisteddef]{Definition}
\begin{chain}\label{chain}
\begin{enumerate}
\item A map $j : \Gamma \to Z(A) \cap A^*$ is a 1-cocycle if 
\begin{equation} \label{eq:jmap}
j(\gamma \mu) = \big ( (j\gamma) \cdot \mu \big ) j \mu  , \quad \,\,\, \forall \gamma, \mu \in \Gamma .
\end{equation}
\item Given a map $j : \Gamma \to A$, a linear functional $\tau :A \to \mathbb{C}$ is said to have the change of variable property with respect to $j$, if
\[ \tau((a \cdot \gamma)j \gamma)= \tau (a),  \quad \,\,\, \forall a \in A, \,\,\, \gamma \in \Gamma.\]
\end{enumerate}
\end{chain}
Notice that \eqref{eq:jmap} amounts to saying that $j$ is a
(multiplicative) group $1-$cocycle for $H^1(\Gamma , Z(A) \cap A^*).$

\newtheorem{automor}[twisteddef]{Proposition}
\begin{automor}\label{twistedderivation}
Let $A$ be an algebra with a right action of a group $\Gamma$ by automorphisms, and $j : \Gamma \to Z(A) \cap A^*$ be a 1-cocycle.
\begin{enumerate}
\item The map $\sigma : A \rtimes \Gamma \to A \rtimes \Gamma$  given by
\[ \sigma (a \otimes \gamma) = \big ( (j \gamma^{-1})\cdot \gamma \big ) a \otimes \gamma\]
is an automorphism.
\item Let $\delta : A \to A$ be a derivation such that $\delta (a \cdot \gamma) = (\delta (a) \cdot \gamma) \, j \gamma$ for all $a \in A, \gamma \in \Gamma$. Then for any $s=1,2, \dots$, the map $\delta'_s : A \rtimes \Gamma \to A \rtimes \Gamma$ defined by
\[ \delta'_s (a \otimes \gamma) = \Big ( \delta \Big ( (a \cdot \gamma^{-1}) (j \gamma^{-1})^s \Big ) (j \gamma^{-1})^{-s} \otimes 1 \Big ) \Big ( 1 \otimes \gamma \Big ) \]
is a $\sigma$-derivation on $A \rtimes \Gamma$. Also $\delta'_s \circ \sigma = \sigma \circ \delta'_s$ if $\delta \circ j = 0$.
\item If $\tau : A \to \mathbb{C}$ is a trace such that $\tau \circ \delta = 0$, then $\tau' \circ \delta'_s =0$ where the linear functional $\tau' : A \rtimes \Gamma \to \mathbb{C}$ is defined by
\[
\tau' (a \otimes \gamma) = 0 \,\, if \,\,  \gamma \neq 1, \,\,\, and \,\,\, \tau' (a \otimes 1) = \tau(a).\]
Also, if $\tau$ has the change of variable property with respect to $j$, then $\tau'$ is a $\sigma$-trace on $A \rtimes \Gamma.$

\end{enumerate}

\begin{proof}
\begin{enumerate}
\item For any $a,b \in A$, $\gamma, \mu \in \Gamma$, we have
\begin{eqnarray}
 \sigma(a \otimes \gamma) \sigma (b \otimes \mu) &=& \big ( (j\gamma^{-1} \cdot \gamma) a \otimes \gamma \big ) \big ( (j\mu^{-1} \cdot \mu) b \otimes \mu \big ) \nonumber \\
 &=& \big (j\gamma^{-1} \cdot (\gamma \mu) \big )  (a \cdot \mu) (j\mu^{-1} \cdot \mu) b \otimes \gamma \mu  , \nonumber
\end{eqnarray}
and
\begin{eqnarray}
\sigma \big ( (a \otimes \gamma)(b \otimes \mu) \big ) &=& \sigma \big ( (a \cdot \mu) b \otimes \gamma \mu \big ) \nonumber \\
&=& \big ( j(\mu^{-1} \gamma^{-1}) \cdot (\gamma \mu) \big ) (a \cdot \mu) b \otimes \gamma \mu \nonumber \\
&=& \bigg ( \Big ( \big ( j(\mu^{-1}) \cdot \gamma^{-1}) j \gamma^{-1} \Big ) \cdot (\gamma \mu) \bigg ) \,\,  (a \cdot \mu) \,\, b \, \otimes \gamma \mu \nonumber \\
&=& (j \mu^{-1} \cdot \mu )\,\, \big ( (j \gamma^{-1})\cdot (\gamma \mu) \big ) \,\, (a \cdot \mu) \, b \otimes \gamma \mu  . \nonumber
\end{eqnarray}
Therefore $\sigma$ is an algebra homomorphism since $j(\Gamma) \subset Z(A)$. It is also easy to see that $\sigma$ is an automorphism since $j(\Gamma) \subset A^*$.

\item For any $a,b \in A$, and $\gamma , \mu \in \Gamma$, we have
\begin{eqnarray} 
&&  \delta'_s((a \otimes \gamma)(b \otimes \mu)) = \delta'_s ((a \cdot \mu)b \otimes \gamma \mu)  \nonumber \\                                             &=& \delta \Big ( \big ( (a \cdot \mu) b \big ) \cdot (\mu^{-1} \gamma^{-1}) \, (j\mu^{-1} \gamma^{-1})^s \Big ) \, (j\mu^{-1} \gamma^{-1})^{-s} \otimes 1 \,\, 1\otimes \gamma \mu \nonumber \\
&=& \delta \Big ( (a \cdot \gamma^{-1}) \, b \cdot (\mu^{-1} \gamma^{-1}) \, (j\mu^{-1} \gamma^{-1})^s \Big ) \, (j\mu^{-1} \gamma^{-1})^{-s} \otimes 1 \,\, 1\otimes \gamma \mu \nonumber \\
&=& \delta(a \cdot \gamma^{-1}) \, b \cdot (\mu^{-1} \gamma^{-1}) \otimes 1 \,\, 1\otimes \gamma \mu \, +\nonumber \\
&& (a \cdot \gamma^{-1}) \, \delta (b \cdot (\mu^{-1} \gamma^{-1})) \otimes 1 \,\, 1\otimes \gamma \mu \, +\nonumber \\
&& (a \cdot \gamma^{-1}) (b \cdot (\mu^{-1} \gamma^{-1})) \delta \big ( (j\mu^{-1} \gamma^{-1})^s \big ) (j\mu^{-1} \gamma^{-1})^{-s} \otimes 1 \,\, 1\otimes \gamma \mu. \nonumber
\end{eqnarray}
Since $(j\mu^{-1} \gamma^{-1})^s = (j\mu^{-1})^s \cdot \gamma^{-1} \, (j\gamma^{-1})^s$, the latter is equal to:
\begin{eqnarray} \label{twistedderivationl}
&& \delta(a \cdot \gamma^{-1}) \cdot (\gamma \mu) \, b \otimes \gamma \mu \, +\nonumber \\
&& (a \cdot \mu) \, \delta (b \cdot (\mu^{-1} \gamma^{-1})) \cdot (\gamma \mu) \otimes \gamma \mu \, +\nonumber \\
&& (a \cdot \mu)\, b\,  \delta ((j\mu^{-1})^s) \cdot \mu \, (j\gamma^{-1})^{s+1} \cdot (\gamma \mu ) \, (j\mu^{-1} \gamma^{-1})^{-s} \cdot (\gamma \mu) \otimes \gamma \mu \, + \nonumber \\
&& (a \cdot \mu)\, b\, (j\mu^{-1})^s \cdot \mu \,\, \delta ((j\gamma^{-1})^{s}) \cdot (\gamma \mu ) \, (j\mu^{-1} \gamma^{-1})^{-s} \cdot (\gamma \mu) \otimes \gamma \mu. 
\end{eqnarray}

Also we have
\begin{eqnarray} \label{twistedderivationr}
&& \delta'_s (a \otimes \gamma)(b \otimes \mu) + \sigma (a \otimes \gamma) \, \delta'_s(b \otimes \mu) \nonumber \\
&=& \Big ( \delta \big ( (a \cdot \gamma^{-1}) (j \gamma^{-1})^s \big ) (j \gamma^{-1})^{-s} \otimes 1 \Big ) \Big ( b \otimes \gamma \mu \Big ) + \nonumber \\
&& \Big ( (j \gamma^{-1} \cdot \gamma) a \otimes \gamma \Big ) \Big ( \delta \big ( (b \cdot \mu^{-1}) (j \mu^{-1})^s  \big ) (j \mu^{-1})^{-s} \otimes 1 \Big ) \Big ( 1 \otimes \mu \nonumber \Big ) \\
&=& \delta (a \cdot \gamma^{-1}) \cdot (\gamma \mu) \, b \otimes \gamma \mu \, + \nonumber \\
&& (a \cdot \mu) \, \delta ((j\gamma^{-1})^s) \cdot (\gamma \mu) \, (j \gamma^{-1})^{-s} \cdot (\gamma \mu) \, b \otimes \gamma \mu \, + \nonumber \\
&& (j \gamma^{-1}) \cdot (\gamma \mu) \, (a \cdot \mu) \, \delta (b \cdot \mu^{-1}) \cdot \mu \otimes \gamma \mu \, + \nonumber \\
&& (j \gamma^{-1}) \cdot (\gamma \mu) \, (a \cdot \mu) \, b \, \delta ((j \mu^{-1})^s) \cdot \mu \, (j \mu^{-1})^{-s} \cdot \mu \otimes \gamma \mu.
\end{eqnarray}
 It is easy to see that $\delta (a), \sigma (a) \in Z(A)$ for any $a \in Z(A)$.
 Therefore by \eqref{twistedderivationl} and \eqref{twistedderivationr},
 $\delta'_s$ is a $\sigma$-derivation. \\

Now assume that $\delta \circ j = 0$. Then we have
\begin{eqnarray}
\sigma \circ \delta'_s (a \otimes \gamma) &=& \sigma \Big (   \delta \big ( (a \cdot \gamma^{-1}) (j \gamma^{-1})^s \big ) \cdot \gamma \, (j \gamma^{-1})^{-s} \cdot \gamma \otimes \gamma   \Big ) \nonumber \\
&=& (j \gamma^{-1} \cdot \gamma) \, \delta \big ( (a \cdot \gamma^{-1}) (j \gamma^{-1})^s \big ) \cdot \gamma \, (j \gamma^{-1})^{-s} \cdot \gamma \otimes \gamma  \nonumber \\
&=& (j \gamma^{-1} \cdot \gamma )  \, \delta (a \cdot \gamma^{-1}) \cdot \gamma \otimes \gamma, \nonumber
\end{eqnarray}
and
\begin{eqnarray}
\delta'_s \circ \sigma (a \otimes \gamma ) &=& \delta'_s \big (  (j \gamma^{-1} \cdot \gamma) a \otimes \gamma  \big ) \nonumber \\
&=& \delta (j \gamma^{-1} \, a \cdot \gamma^{-1} \,  (j \gamma^{-1})^s ) \cdot \gamma \, (j \gamma^{-1})^{-s} \cdot \gamma \otimes \gamma \nonumber \\
&=& \big ( j \gamma^{-1} \, \delta (a \cdot \gamma^{-1})  \big ) \cdot \gamma \otimes \gamma  .\nonumber
\end{eqnarray}
So $\sigma$ and $\delta'_s$ commute if $\delta \circ j = 0$.

\item It is easy to see that $j1 =1$, from which it follows that $\tau ' \circ \delta'_s = 0$ if $\tau \circ \delta = 0$.

To prove that $\tau'$ is a $\sigma$-trace, we have
\[ (a \otimes \gamma)(b \otimes \mu) = (a \cdot \mu) b \otimes \gamma \mu, \]
and
\begin{eqnarray}
\sigma (b \otimes \mu) \,\, (a \otimes \gamma) &=& \big ( (j\mu^{-1} \cdot \mu) b \otimes \mu \big ) (a \otimes \gamma ) \nonumber \\
&=& \big ( j\mu^{-1} \cdot (\mu \gamma) \big ) (b \cdot \gamma) a \otimes \mu \gamma. \nonumber
\end{eqnarray}
Therefore if $\mu \gamma \neq  1$, then $\tau'$ is $0$ on both of the above terms, and if $\mu \gamma =  1$ then
\[ \tau'((a \otimes \gamma)(b \otimes \mu) ) = \tau ((a\cdot \mu)b) \]
and
\begin{eqnarray}
\tau' ( \sigma (b \otimes \mu) \,\, (a \otimes \gamma)) &=& \tau ((j \mu^{-1})(b \cdot \mu^{-1})  a) \nonumber \\
 &=& \tau \big ( (j \mu^{-1}) (  ( b  (a \cdot \mu) ) \cdot \mu^{-1} \big ) \nonumber \\
 &=& \tau (b  (a \cdot \mu) ). \nonumber
\end{eqnarray}
Therefore $\tau'$ is a $\sigma$-trace on $A \rtimes \Gamma$.
\end{enumerate}
\end{proof}
\end{automor}

\begin{example}\rm{
\begin{enumerate}
\item Let $M$ be a smooth oriented manifold and $\omega$ a volume form on $M$. Let $\Gamma = Diff(M)$ be the group of diffeomorphisms of $M$. The map $j: \Gamma \to C^{\infty} (M)$ defined by
\[ \gamma ^* (\omega) = j(\gamma) \omega \]
is easily seen to be a 1-cocycle.

\item Let $\chi : \Gamma \to \mathbb{C}^*$ be a 1-dimensional character of a group $\Gamma$ which acts on an algebra $A$ by automorphisms, and let $j(\gamma) = \chi (\gamma) 1_{A}$. Then $j$ is a 1-cocycle,
and a derivation $\delta : A \to A$ is compatible with $j$ if and only if
$\delta (a \cdot \gamma) = \chi (\gamma) \,\, \delta(a) \cdot \gamma$ for
any $a \in A$, and $\gamma \in \Gamma.$

\end{enumerate}}
\end{example}

\newtheorem{connectiondef}[twisteddef]{Definition}
\begin{connectiondef}\label{connectiondef}
Let $A$ be an algebra, $\delta : A \to A$ a derivation, and $E$ a left $A$-module. A linear map $\nabla : E \to E$ is said to be a connection if it satisfies the Leibniz rule, i.e.

\[ \nabla (a \, \xi) = \delta (a) \, \xi + a \, \nabla (\xi), \,\,\, \forall a \in A, \,\,\, \xi \in E. \]

If $\sigma : A \to A$ is an automorphism and $\delta : A \to A$  a $\sigma$-derivation, then a linear map $\nabla : E \to E$ is said to be a twisted connection if it satisfies the twisted Leibniz rule:
\[ \nabla (a \, \xi) = \delta (a) \, \xi + \sigma (a) \, \nabla (\xi), \,\,\, \forall a \in A, \,\,\, \xi \in E. \]
\end{connectiondef}
The  notion of twisted connection was used by  Polishchuk \cite{pol} in his study of noncommutative
toroidal orbifolds. \\

One of the conditions for a twisted spectral triple $(A, H, D)$ in the sense of
Connes-Moscovici \cite{conmos} is the boundedness of the 
operators $[D, a]_{\sigma} := Da - \sigma(a) D$. Given a twisted connection $\nabla$, one can try to define a twisted  spectral triple by letting $D = \nabla$. Then $\nabla a - \sigma (a) \nabla = \delta (a)$ shows the boundedness condition 
is satisfied provided $\delta (a)$ acts by a bounded operator. \\

In the following proposition, $A$ is an algebra endowed with a right action of a
group $\Gamma$ by automorphisms, with representations $\pi : A \to \text{End} \,(E)$, and
$\rho: \Gamma \to \text{GL} (E)$,  defining  a covariant system, i.e.
\[\pi (a \cdot \gamma ) = \rho (\gamma^{-1}) \, \pi(a) \, \rho(\gamma), \,\,\, \forall a \in A, \,\,\, \gamma \in \Gamma.  \]
Then  we obtain  a representation $\pi' : A \rtimes \Gamma \to \text{End} (E)$ given by
\[\pi' (a \otimes \gamma) = \rho(\gamma) \pi (a), \,\,\,  \forall a \in A, \,\,\, \gamma \in \Gamma.\]
Also let  $j : \Gamma \to Z(A) \cap A^*$ be a 1-cocycle and 
$\delta : A \to A$  a derivation  as in Proposition \ref{twistedderivation}.
Therefore we have an automorphism $\sigma : A \rtimes \Gamma \to A \rtimes
\Gamma $ and we fix a $\sigma$-derivation $\delta'_s : A \rtimes \Gamma \to A \rtimes \Gamma$ for some $s \in \mathbb{N}$.

\newtheorem{twistedconnection}[twisteddef]{Proposition}
\begin{twistedconnection}\label{twistedconnection}
A connection  $\nabla: E \to E $ for $A$ is a twisted connection for $A \rtimes \Gamma$ if and only if
\begin{equation} \label{twistedconncond}
\nabla  \rho (\gamma) = \rho (\gamma) \Big ( (s-1) \, \pi \big (\delta (j \gamma^{-1}\cdot \gamma) \big ) + \nabla \pi (j \gamma^{-1}\cdot \gamma)  \Big ), \,\,\, \forall \gamma \in \Gamma.
\end{equation}

\begin{proof}
Since $\sigma(a \otimes 1) = a \otimes 1$, and $\delta'_s (a \otimes 1) = \delta(a) \otimes 1$, it suffices to show that
\[ \nabla \, \pi'(1 \otimes \gamma)    = \pi' (\delta'_s (1 \otimes \gamma))  + \pi' (\sigma(1\otimes \gamma)) \nabla \]
for all $\gamma \in \Gamma$, if and only if \eqref{twistedconncond} holds.
First we compute $\pi' (\delta'_s (1 \otimes \gamma))$:
\begin{eqnarray}
\pi' (\delta'_s (1 \otimes \gamma)) &=& \pi' \Big (   \delta \big ((j\gamma^{-1})^s \big ) \cdot \gamma \, (j\gamma^{-1})^{-s} \cdot \gamma \otimes \gamma \Big ) \nonumber \\ 
&=& \rho (\gamma) \, \pi \Big (  \delta \big ((j\gamma^{-1})^s \big ) \cdot \gamma \, (j \gamma)^s \Big ) \nonumber \\
&=& \rho (\gamma) \, \pi \Big (  \delta \big ((j\gamma^{-1})^s \big ) \cdot \gamma \, j \gamma \, (j \gamma)^{s-1}\Big ) \nonumber \\
&=& \rho (\gamma) \, \pi \Big (  \delta \big ((j\gamma^{-1})^s \cdot \gamma \big )  \, (j \gamma)^{s-1} \Big ) \nonumber \\
&=& \rho (\gamma) \, \pi \Big ( s \, (j\gamma^{-1} \cdot \gamma)^{s-1} \delta (j\gamma^{-1} \cdot \gamma ) (j\gamma)^{s-1} \Big ) \nonumber \\
&=& s \, \rho (\gamma) \, \pi \big ( \delta (j \gamma^{-1} \cdot \gamma) \big ). \nonumber
\end{eqnarray}
Now we have 
\begin{eqnarray}
&& \pi' (\delta'_s (1 \otimes \gamma))  + \pi' (\sigma(1\otimes \gamma)) \nabla \nonumber \\
&=& (s-1) \, \rho (\gamma) \, \pi \big ( \delta (j \gamma^{-1} \cdot \gamma) \big ) + \rho (\gamma) \Big ( \pi \big ( \delta (j \gamma^{-1} \cdot \gamma ) \big ) + \pi (j \gamma^{-1} \cdot \gamma) \, \nabla \Big ) \nonumber \\
&=&\rho (\gamma) \Big ( (s-1) \, \pi \big (\delta (j \gamma^{-1}\cdot \gamma) \big ) + \nabla \pi (j \gamma^{-1}\cdot \gamma)  \Big ). \nonumber
\end{eqnarray}

Therefore $\nabla$ is a twisted connection if and only if \eqref{twistedconncond} holds.

\end{proof}

\end{twistedconnection}

\begin{remark}\rm{
Let $s > 0$ be a real number and assume $\big ( j (\gamma) \big )^s \in A$ is defined for all $\gamma \in \Gamma$. Propositions \ref{twistedderivation} and \ref{twistedconnection} continue to hold for these values of $s$ as well. For this we need the extra condition $\delta (x^s) = s x^{s-1} \delta (x)$ to hold for all $x = j(\gamma) $, $\gamma \in \Gamma$.}
\end{remark}

\begin{example}\rm{ Let $ C^{\infty} (S^1)$ be the algebra of smooth functions on the circle $S^1 = \mathbb{R} / \mathbb{Z}$, and $\Gamma \subset Diff^+ (S^1)$ a group of orientation preserving diffeomorphisms of the circle as in \cite{conmos}.  We represent the algebra $C^{\infty} (S^1)$ by bounded operators in the Hilbert space $L^2(S^1)$ by
\[ (\pi (g) \, \xi) (x) = g(x) \, \xi (x), \,\,\, \forall g \in C^{\infty} (S^1), \,\,\,  \xi \in L^2(S^1), \,\,\, x \in \mathbb{R} / \mathbb{Z}. \]
Define a representation of $\Gamma$ by bounded operators in  $L^2(S^1)$ by
\[ (\rho(\phi^{-1}) \, \xi )(x) = \phi'(x)^{\frac{1}{2}} \, \xi(\phi(x)), \,\,\, \forall \phi \in \Gamma, \,\,\,  \xi \in L^2(S^1), \,\,\, x \in \mathbb{R} / \mathbb{Z}. \]
The group $\Gamma$ acts on $ C^{\infty} (S^1)$ from right by composition and one can
easily check  that the above representations give a covariant system which yields the
representation of $ C^{\infty} (S^1) \rtimes \Gamma$ as  in  \cite{conmos}.
The map $j : \Gamma \to C^{\infty} (S^1)$ defined by $j(\phi) = \phi'$ is a 1-cocycle and the derivation $\delta: C^{\infty} (S^1)
\to C^{\infty} (S^1), \delta (f) = \frac{1}{i} f'$ is compatible with $j$. Now by using
Proposition \ref{twistedderivation}, we obtain an automorphism $\sigma$
of $ C^{\infty} (S^1) \rtimes \Gamma$ which agrees with the automorphism
in \cite{conmos}, and a twisted derivation $\delta'_{\frac{1}{2}}$. Note that since it is possible to take the square root of the elements in the image of $j$ in this example, we can let $s=\frac{1}{2}$ to obtain a twisted derivation. Now if we let $\nabla = \frac{1}{i} \frac{d}{dx}$, one can see that the equality
\eqref{twistedconncond} holds, therefore $\nabla$ is a twisted connection for $C^{\infty}(S^1) \rtimes \Gamma$  by Proposition \ref{twistedconnection}.

Also the linear map $\tau : C^{\infty}(S^1) \to \mathbb{C}$ defined by
\[ \tau (g) = \int_{\mathbb{R} / \mathbb{Z}} g(x) \, dx, \,\,\, \forall g \in C^{\infty}(S^1), \]
is a trace which has the change of variable property with respect to $j$, and $\tau \circ \delta = 0$. Therefore by Proposition \ref{twistedderivation}, one obtains a twisted trace $\tau' : C^{\infty}(S^1) \rtimes \Gamma \to \mathbb{C} $ such that $\tau' \circ \delta'_{\frac{1}{2}} = 0.$

}
\end{example}

\end{document}